\documentclass{article}
\usepackage{graphicx} 
\usepackage[a4paper,top=3cm,bottom=2cm,left=3cm,right=3cm,marginparwidth=1.75cm]{geometry}
\usepackage{amsmath}
\usepackage[title]{appendix}
\usepackage[T1]{fontenc}
\usepackage{yfonts}
\usepackage{enumitem,xparse}
\usepackage{amsmath}
\usepackage{amssymb}
\usepackage{hyperref}
\hypersetup{colorlinks=true, citecolor=green}
\usepackage{alltt}
\usepackage{setspace}
\usepackage[english]{babel}
\usepackage{euscript}
\usepackage{yfonts}
\usepackage{tikz}
\usepackage{amstext}
\usepackage{euscript}
\usepackage{amsfonts}
\usepackage{dcolumn}
\usepackage{bm}
\usepackage{psfrag}
\usepackage{euscript}
\usepackage{array}
\usepackage{multirow}
\usepackage{supertabular}
\newcommand{\old}[1]{{}}
\usepackage{theoremref}
\usepackage{amsthm}
\usepackage[english]{babel}
\usepackage{textcomp}
\usepackage{amssymb}
\usepackage{amsmath}
\usepackage{authblk}
\usepackage{tipa}
\usepackage{graphicx}
\usepackage{enumerate}
\usepackage{dsfont}
\usepackage{mathrsfs} 
 \usepackage{amsmath}

\theoremstyle{definition}
\newtheorem{definition}{Definition}[section]
\newtheorem{remark}{Remark}[section]
 \usepackage{xcolor}
\newtheorem{theorem}{Theorem}[section]
\theoremstyle{lemma}

\theoremstyle{proposition}
\newtheorem{proposition}[theorem]{Proposition}
\newtheorem{corollary}{Corollary}[section]

\newcommand\numberthis{\addtocounter{equation}{1}\tag{\theequation}} 
\let\svvert\|
\def\norm#1{\left\svvert#1\right\svvert}
\def\|#1\|{\norm{#1}}

\title{Application of Onsager's Variational Principle to the Navier-Stokes Equations}

\date{}

\author{Hamid A. Said\footnote{Department of Mathematics, Kuwait University, Safat 13060, Kuwait.  $\mathtt{ hamids@sci.kuniv.edu.kw}.$}}

\begin{document}

\maketitle
\begin{abstract}
In this note we propose a basic $L^2$-based approach for studying  the global and local energy equalities of the incompressible 3D Navier-Stokes equations in the standard energy class on $\mathbb{T}^3 \times (0,T]$. Motivated by L. Onsager's principle of least dissipation of energy \cite{onsager1931reciprocal1, onsager1931reciprocal2}, we give a new sufficient condition for the energy equalities in terms of the limit of a sequence of minimizers. In particular, we show that the equalities are attained if this limit is non-vanishing. We observe that the indeterminacy in the vanishing case is reminiscent of turbulence-driven energy transfer dominating other transport processes.
 \end{abstract}

\section{Introduction} \label{introduction}
Consider the incompressible 3D Navier-Stokes equations
\begin{equation}\label{NS}
  \left\{
  \begin{array}{ll}
  \vspace*{0.1 in}
   \partial_t u^\nu + (u^\nu \cdot \nabla) u^\nu  + \nabla p^\nu - \nu \Delta u^\nu = 0 \, ,\\
    \nabla \cdot u^\nu = 0\, ,
    \end{array}
     \right.
\end{equation}
posed over $\mathbb{T}^3 \times (0,T]$. The constant $\nu >0$ is the kinematic viscosity, and the velocity-pressure pair $(u^\nu, p^\nu)$ is assumed to have periodic boundary conditions with zero mean over $x\in \mathbb{T}^3$. Unless otherwise stated, we will consider the Navier-Stokes equations for fixed $\nu$, and so we set $u \doteq u^\nu$. 
\smallskip
\\
\indent In his seminal work \cite{leray1934mouvement}, J. Leray proved the existence of at least one global in-time (\emph{i.e.} for all $T >0$) weak solution to \eqref{NS} in $\mathbb{R}^3$. Approximately 17 years later, E. Hopf \cite{hopf1950anfangswertaufgabe} proved a similar result over smooth bounded domains with Dirichlet boundary conditions. Both solutions were shown, in addition, to satisfy the \emph{global energy inequality}:

\begin{equation} \label{NS_EN_ineq}
    \dfrac{1}{2}\lVert u(t) \rVert^2 + \nu \int_0^t \lVert \nabla u \rVert^2 ds \leq \dfrac{1}{2}\lVert u_0 \rVert^2 \, ,
\end{equation}
for all $t \geq 0$, where $\lVert \cdot \rVert$ denotes the $L^2$-norm. Solutions satisfying the above energy inequality have come to be known as \emph{Leray-Hopf solutions}, and have been a subject of extensive study for many decades. We refer the reader to the classical works \cite{constantin1988navier, temam2024navier} and the more recent texts \cite{robinson2016three, boyer2012mathematical} for comprehensive treatments on the subject.
\medskip
\\
\indent The other energy inequality of interest is the \emph{local energy inequality}, established first by V. Scheffer \cite{scheffer1977hausdorff} for Leray solutions:
\begin{align} \label{local-en-ineq-standard-intro} 
 2 \nu \int \int | \nabla u|^2 \phi\, dx dt \leq  \int \int | u|^2 \left( \partial_t \phi + \nu \Delta \phi \right) \, dx dt + \int \int  \left( | u|^2 + 2p \right)(u \cdot \nabla) \phi \, dx dt \, ,
\end{align}
which holds for any non-negative $C^\infty$ scalar function $\phi$ compactly supported on $\mathbb{R}^3 \times \mathbb{R}$. We can rewrite the local energy inequality \eqref{local-en-ineq-standard-intro} in more convenient form:
\begin{align} \label{local-en-ineq-intro}
       \partial_t \left( \frac{|u|^2}{2}  \right) + \nabla \cdot \left(u \left( \frac{|u|^2}{2}  + p\right) - \nu \nabla \left( \frac{|u|^2}{2} \right) \right)  \leq  - \nu (\nabla u)^2 \, ,
\end{align}
where the inequality is understood to hold in the sense of distributions\footnote{The result remains valid in the case of bounded domains \cite{caffarelli1982partial} and  the absence of boundaries \cite{Ibrahim007construction}. It remains unknown whether Hopf solutions, obtained through a Galerkin scheme, satisfy the local energy inequality. }. The inequality is well-defined, since the regularity class of weak solutions $u$ is $ L^p(0,T; L^q)$, \, $\frac{2}{p} + \frac{3}{q} = \frac{3}{2}, \, \, 2\leq q \leq 6$, and the scalar pressure field $p(x,t)$ belongs to $L^{\frac{3}{2}}(0,T;L^{\frac{3}{2}})$. We note, excursively, that the above  energy inequality plays a central role in the pioneering works of Scheffer \cite{scheffer1977hausdorff} and, later, of Caffarelli, Kohn, and Nirenberg \cite{caffarelli1982partial} on the partial regularity of the weak solutions to the Navier-Stokes equations -- known as \emph{suitable} weak solutions. 
\medskip
\\
\indent A particular area of inquiry related to Leray-Hopf solutions centers around the energy equalities, namely the \emph{global energy equality}
\begin{equation} \label{global-en-eq-intro}
    \dfrac{1}{2}\lVert u(t) \rVert^2 + \nu \int_0^t \lVert \nabla u \rVert^2 ds = \dfrac{1}{2}\lVert u_0 \rVert^2 \, ,
\end{equation}
and the \emph{local energy equality}
\begin{align} \label{local-en-eq-intro}
       \partial_t \left( \frac{|u|^2}{2}  \right) + \nabla \cdot \left(u \left( \frac{|u|^2}{2}  + p\right) - \nu \nabla \left( \frac{|u|^2}{2} \right) \right)  =  - \nu (\nabla u)^2 \, .
\end{align}
Specifically, one seeks conditions, or ideally a sharp regularity threshold, under which one or the other of the above two equations hold. Various sufficient conditions have been established in the literature. Earliest among them was found by J. L. Lions \cite{lions1960regularite} who proved that if weak solutions, in addition, belonged to $L^4(0,T; L^4)$, then the energy equality \eqref{global-en-eq-intro} would hold. Shinbrot \cite{shinbrot1974energy} on the other hand obtained a more general result through interpolation: the global energy equality is valid for solutions $u \in L^p(0,T; L^q)$ with $\frac{1}{p} + \frac{1}{q} \leq \frac{1}{2}$ and $q \geq 4$. It is worth noting that Shinbrot's condition is not scale invariant\footnote{Compare for instance with the scaling given by the Serrin condition (also termed Ladyzhenskaya-Prodi-Serrin conditions): if $u \in L^p(0,T; L^q)$ with $\frac{2}{p} + \frac{3}{q} \leq 1$ and $q \geq 3$, then $u$ is smooth over the interval $(0,T]$. Clearly, $\frac{2}{p} + \frac{2}{q} \leq \frac{2}{p} + \frac{3}{q} \leq 1$.  }. A weaker condition for energy conservation follows from the result of Cheskidov, Constantin, Friedlander and Shvydkoy \cite{cheskidov2008energy} in  an $L^3$-based Besov space: if a weak solution $u \in L^3 (0,T; B^{\frac{1}{3}}_{3,\infty})$, then equation \eqref{local-en-eq-intro} holds\footnote{The norm in this Besov space is given by $\lVert v \rVert_{B^\beta_{p,\infty}} = \lVert v \rVert_{L^p} + \sup_{|y| > 0} \dfrac{\lVert v(\cdot + y) - v(\cdot) \rVert_{L^p}}{|y|^\beta}$ .} (see also \cite{cheskidov2020energy}). Combining this condition with the commutator estimate of Constantin, E, and Titi \cite{constantin1994onsager} (originally proven for the incompressible Euler equations; that is for $\nu := 0$ in equations \eqref{NS}), the global and local energy equations for Leray solutions are proved to hold in $L^3 (0,T; B^{\beta}_{3,\infty})$ for $\beta \geq \frac{1}{3}$. We also refer to \cite{berselli2023energy, cheskidov2020energy, leslie2018conditions, shvydkoy2009geometric} for other and more recent partial results in this area. As remarked in \cite{shvydkoy2009geometric}, for an appropriately chosen sequence of test functions, the global energy equality \eqref{global-en-eq-intro} can be recovered from equation \eqref{local-en-eq-intro} for any weak solution.
\medskip
\\
\indent For fixed values of viscosity,  J. Duchon and R. Robert \cite{duchon2000inertial} successfully quantified the dissipation of energy that could occur due to potential singularities in Leray-Hopf solutions. In fact, it was proven that any weak solution $u$ to the Navier-Stokes equations satisfies the following local energy identity point-wise in $(x,t) \in \mathbb{T}^3 \times [0,T]$ (in the sense of distributions):
\begin{align} \label{local_en_D_intro}
    \partial_t \left( \frac{|u|^2}{2}  \right) + \nabla \cdot \left(u \left( \frac{|u|^2}{2}  + p\right) - \nu \nabla \left( \frac{|u|^2}{2} \right) \right)  =  - \nu (\nabla u)^2 - D(u) \, ,
\end{align}
where the distribution $D(u)(x,t)$ is given by the relation
\begin{align*} 
    D(u)(x,t) = \lim_{n \to \infty} \dfrac{1}{4} \int_{\mathbb{T}^3} \nabla \eta^{n} (y) \cdot \left( u (x+y,t) - u(x,t)\right) |u (x+y,t) - u (x,t)|^2 dy \, ,
\end{align*}
and function $\eta$ is a non-negative, radial, smooth function with compact support and unit mass on $\mathbb{R}^3$, and $\eta^n (y) = \frac{1}{\delta_n^3} \eta\left(\frac{y}{\delta_n}\right)$ with  $\delta_n \to 0$ as $n \to \infty$. The above limit is shown to hold regardless of the choice of mollifier $\eta$\footnote{Clearly, for sufficiently regular solutions $u$ the limit: $\lim_{|y| \to 0} \int_{\mathbb{T}^3} |u (x+y,t) - u (x,t)|^3 dy$ vanishes leading to $D(u) = 0$.}. In terms of mollified velocity $\overline{u}^{\,n}$, the limit can be re-written as
$$ D(u) = \lim_{n \to \infty} (- \nabla  \overline{u}^{\,n} :  \overline{R}^{\, n}) \, ,$$
which holds in the sense of distributions\footnote{ Observe that $\lim_n [\overline{u}^{ \, n} \cdot (\nabla\cdot \overline{R}^{ \, n})] = \lim_n (- \nabla \overline{u}^{ \, n} : \overline{R}^{ \, n})$ in the sense of distributions, and thus we get $\lim_n \partial_i ( \overline{u}^{\,n}_j \,  \overline{R}^{\, n}_{ij}) = 0$ in the distributional sense.}, and where we have defined $\overline{R}^{ \, n} = (\overline{u \otimes u})^{n} - \overline{u}^{ \, n} \otimes \overline{u}^{ \, n}$.  As alluded to earlier, in the case where $u$ is taken to be a Leray solution we have $D(u) \geq 0$. Additionally, one can readily verify that establishing the global and local energy equalities are equivalent for Leray solutions. If $D(u)$ vanishes identically for Leray-Hopf solutions, then this signifies that viscosity alone is the driving mechanism behind the dissipation of energy. Even beyond the question of identifying the dissipative mechanism of the energy of the flow, establishing the energy equalities play a key role in connecting statistical solutions of the Navier-Stokes equations -- first introduced by C. Foias \cite{foiacs1972statistical, foias2001navier} -- to the Kolmogorov $(\frac{3}{5})$ spectrum in the theory of turbulence. These solutions consist of a family of time-dependent, spatially homogeneous probability measures $\mu_t$ (for $t \in [0,\infty)$) on an appropriate Hilbert space, satisfying the Navier-Stokes equations in an averaged sense. If the \emph{averaged} global energy equality is then satisfied, that is equation \eqref{global-en-eq-intro} when averaged against the probability measure $\mu_t$, then the total averaged kinetic energy displays universal features characteristic of homogeneous isotropic turbulence: long-time decay independent of the initial state of the flow, and an energy spectrum in Fourier space that follows the Kolmogorov power law $\kappa^{-\frac{3}{5}}$, where $\kappa$ denotes the wavenumber. See also \cite{CHEN2024134274} for a connection to the Kolmogorov scaling  when the energy equation is beyond reach. This particular approach is based on an entropy maximum principle formulated over an energy-enstrophy surface, defined in \eqref{surface}, with a prescribed energy profile $e(t)$. Further discussion can be found in Section \ref{discussion}.
\medskip
\\
\indent  While it is clear that, in principle, further regularity of weak solutions is not needed for the energy equalities to hold, it has remained an open question whether Leray-Hopf solutions do in fact satisfy either the local or global energy equalities.

\subsection{Onsager's principle and statement of main result } \label{main result}
\indent We propose a different approach for studying the energy equalities, which is motivated by L. Onsager's principle of the least dissipation of energy for linear irreversible thermodynamics \cite{onsager1931reciprocal1, onsager1931reciprocal2}.  Onsager argues that physical processes in a thermodynamic system not too far from its equilibrium state must satisfy a maximum principle given in terms of the rate of the internal entropy production function and a dissipation function. The associated Euler-Lagrange equations produce none other than the celebrated Onsager reciprocal relations, and as a consequence give rise to linear equations connecting the thermodynamic forces with the fluxes analogous to equations \eqref{poisson-intro}. Unlike conventional variational principles in mechanics and field theory, such as Hamilton's least action principle, where the evolution equations can be realized as extreme or stationary points of a particular energy functional, Onsager's variational principle selects admissible irreversible thermodynamic processes not too far from the equilibrium state; effectively governing the dissipative macroscopic dynamics in this regime. An outline of the principle of the least dissipation of energy is presented in Appendix \ref{Appendix A}.
\medskip
\\
\indent From a mathematical standpoint, we introduce the following sequence of functionals: 

\begin{align} \label{functional-intro}
    \mathcal{K}^n(v) = \int_0^T \frac{1}{2} \lVert \nabla v(t) \rVert^2 - \langle J^n (t), \nabla v (t) \rangle\, dt \, ,
\end{align}
which, at each fixed $n$, represents the difference between the dissipation function and the entropy production rate for an incompressible Navier-Stokes fluid. See Appendix \ref{Appendix A} for more details. The functions $J^n(x,t) $ express the dissipative fluxes acting within the fluid such as viscosity and possibly another dissipative mechanism driven by the nonlinearity:
\begin{align} \label{sigma}
    J^n(x,t) = \nu \nabla \overline{u}^{\, n}(x,t) - \overline{R}^{ \, n}(x,t)\, .
\end{align}
In this context, one possible mathematical formulation of Onsager's principle is as follows. Let $H$ be the space of divergence free, zero mean periodic vector fields in $L^2(\mathbb{T}^3)$, and let $V$ be its counterpart in $H^1(\mathbb{T}^3)$. Now, consider the following minimization problem (for fixed $n \in \mathbb{N}^+)$: 
\begin{align} \tag{MP}\label{min-prob-intro}
    \inf \lbrace \mathcal{K}^n(v) : v \in \mathcal{B} \rbrace \, ,
\end{align}
where the admissible set is taken to be the \emph{enstrophy ball} $\mathcal{B} = \lbrace v \in L^2(0,T; V) :  \int_0^T \lVert \nabla v \rVert^2dt \leq \frac{1}{2} \lVert u_0 \rVert^2 \rbrace$ for some non-zero $u_0 \in L^2 (\mathbb{T}^3)$. We are now in position to state the main result of the paper
\begin{theorem}  \label{Th1}
\emph{Define $J^n(x,t)$ as in equation \eqref{sigma}, and let  $v^n_*$ be the unique solution  to the problem \eqref{min-prob-intro}. Then there exists a subsequence (not relabeled)  and  $v_* \in L^2(0,T;V)$ such that  $v^n_* \rightharpoonup v_*$ in $L^2(0,T;V)$. If $v_* \neq 0$, then any Leray-Hopf solution $u \in L^\infty(0,T; H) \cap L^2(0,T; V) $ satisfies the local energy equality \eqref{local-en-eq-intro} on $\mathbb{T}^3 \times [0,T)$ in the sense of distributions, hence also the global energy equality \eqref{global-en-eq-intro} for all $t \geq 0$.}
\end{theorem}

\indent The proof of the theorem is detailed in Sections \ref{Variational analysis} and \ref{passage}.  First, however, a few remarks are in order. The application of Onsager's principle in this setting has yielded a convergent (sub)sequence of minimizers  $v^n_*$ for the functional $\mathcal{K}^n$ defined in \eqref{functional-intro}. These functions are related to weak solutions of the Navier-Stokes equations when $v_*$, which is the $L^2(0,T;V)$ weak limit of $ v^n_*$, is non-zero. In particular, under this requirement and by Proposition \ref{prop-w-conv}, $(1-2\lambda_n) v^n_* \rightharpoonup \nu u $ in $L^2(0,T;V)$, where $\lambda_n$ represents a convergent sequence of the Lagrange multipliers associated with the boundary of $\mathcal{B}$; see Section \ref{Minimization problem}. In other words, after passing to the limit, $v^*$ is not identified with $u$ itself, but rather with a constant multiple of $u$, determined by the limit of the sequence of Lagrange multipliers. If instead we have that the miminizers $v^n_*$ lie in the interior of enstrophy ball for all $n$ (or at least for infinitely many $n \in \mathbb{N}$), then $\lambda_n = 0$ and the limit is realized as a genuine weak solution to the equation, \emph{i.e.} $v_* = u$, when $\nu = 1$\footnote{When $\nu \neq 1,$ one can still recover $v_* = u$ simply by introducing the prefactor $\nu$ in front of the quadratic term in functional \eqref{functional-intro}. }. This limitation might be related to data available at the boundary of $\mathcal{B}$, and is further discussed in Section \ref{discussion}. Nevertheless, whether the limit is realized as a weak solution or a constant multiple of it, in either case  the energy identities are still recovered and the sequence of minimizers $v^n_*$ acts as an approximation to said solution; provided that the limit of $v^n_*$ does not vanish. One notable advantage of the variational approach is that it is formulated in the context of the standard energy class $L^\infty(0,T; H) \cap L^2(0,T; V)$, and the requirement $v_* \neq 0$ is well-defined in this class without necessitating additional regularity assumptions about the solution $u$, which is a departure from previous results in this area. Whether the condition $v_* \neq 0$ itself implies higher regularity might be of independent interest. 
\smallskip
\\
\indent We, now, examine the basic ideas behind the proof and then further contextualize our result by making a few additional remarks.

\subsection{Method of proof} \label{proof-sketch}

Let $u$ be a Leray-Hopf solution on $\mathbb{T}^3 \times [0,T]$. The mollification of $u$ admits a sequence of smooth solutions to the \emph{coarse-grained} Navier-Stokes equations (see Section \ref{mollification}):

    \begin{equation}\label{NSCG-intro}
  \left\{
  \begin{array}{ll}
  \vspace*{0.1 in}
   \partial_t \overline{u}^{ \, n} + (\overline{u}^{ \, n} \cdot \nabla) \overline{u}^{ \, n} + \nabla \overline{p}^{ \, n} - \nu \Delta \overline{u}^{ \, n} =  - \nabla \cdot \overline{R}^{ \, n} \, , \\
    \nabla \cdot \overline{u}^{ \, n} = 0\, ,
    \end{array}
     \right.
\end{equation}
where $\overline{R}^{ \, n}(x,t) = (\overline{u \otimes u})^{n} - \overline{u}^{ \, n} \otimes \overline{u}^{ \, n}$. By taking the inner product of equations \eqref{NSCG-intro}$_1$ with $2  \overline{u}^{ \, n} \phi$ (for some scalar test function $\phi$ with support on $\mathbb{T}^3 \times (0,T)$), integrating in space-time, and invoking the divergence free condition, we arrive at the following local energy equation:
\begin{align*} \label{NSCG_LOC_EN} 
 2 \int_0^T \int_{\mathbb{T}^3} \overline{u}^{ \, n} (\nabla \cdot \overline{R}^{ \, n}) \, \phi \, dx dt = \int_0^T \int_{\mathbb{T}^3}  | \overline{u}^{ \, n}|^2 \left( \partial_t \phi + \nu \Delta \phi \right) \, dx dt + &\int_0^T \int_{\mathbb{T}^3}   \left(| \overline{u}^{ \, n}|^2 + 2\overline{p}^{ \, n} \right)(\overline{u}^{ \, n} \cdot \nabla) \phi \, dx dt  \\ 
 &-2 \nu  \int_0^T \int_{\mathbb{T}^3}  | \nabla \overline{u}^{ \, n}|^2 \phi\, dx dt  \, . \numberthis
\end{align*}
Given that the mollification can be removed strongly on the RHS of equation \eqref{NSCG_LOC_EN} (see, for example, \cite{duchon2000inertial}), a plausible approach to establishing the local energy equation \eqref{local-en-eq-intro} therefore is showing that the LHS of equation \eqref{NSCG_LOC_EN}, while admittedly convergent, must vanish in the limit of some subsequence. That is,
\begin{align} \label{NSCG_LOC_EN_cond}
  \int_0^T \int_{\mathbb{T}^3} \overline{u}^{ \, n} ( \nabla \cdot \overline{R}^{ \, n} ) \, \phi \, dx dt \to 0 \, ,
\end{align}
as $n \to \infty$. Since $\overline{u}^{ \, n}$ converges strongly in $L^2(0,T; V)$; then the limit \eqref{NSCG_LOC_EN_cond} can be established if, for instance, one can show that $\nabla \cdot \overline{R}^{ \, n}  \rightharpoonup  0$ in $L^2(0,T; V')$, where $V'$ is the natural dual of $V$.

\begin{remark}
If we, instead, test equations \eqref{NSCG-intro}$_1$ against $\overline{u}^{ \, n}$ only, we obtain the following resolved global energy balance:
\begin{align} \label{NSCG_Glob_EN}
\frac{1}{2} \lVert \overline{u}^{ \, n} (T) \rVert^2 - \frac{1}{2} \lVert \overline{u}^{ \, n} (0) \rVert^2  = - \nu \int_0^T \lVert \nabla \overline{u}^{ \, n}  \rVert^2 dt + \int_0^T \langle \overline{R}^{ \, n} , \nabla \overline{u}^{ \, n}  \rangle \, dt \, ,
\end{align}
where $\langle \cdot, \cdot \rangle$ denote the inner product in $L^2$ space. So similarly for this case as well, if $\nabla \cdot \overline{R}^{ \, n}  \rightharpoonup  0$ in $L^2(0,T; V')$ as $n \to \infty$ (or $\overline{R}^{ \, n}  \rightharpoonup 0$  in $L^2(0,T; H)$  as $n \to \infty$), one can directly prove the global energy equation \eqref{global-en-eq-intro} from the energy balance \eqref{NSCG_Glob_EN}. The $L^2$-in-time weak continuity of Leray-Hopf solutions ensures that the limit can be passed on the LHS of equation \eqref{NSCG_Glob_EN} (see Remark~\ref{C_w remark}).  
\end{remark}

\indent We observe that testing equations \eqref{NSCG-intro}$_1$ with a smooth divergence free vector function $\varphi(x,t)$, and using the strong convergence properties of the mollified velocity on the LHS already gives
\begin{align} \label{lim_R^n}
  \lim_{n \to \infty}  \int_0^T \langle \nabla \cdot \overline{R}^{ \, n} , \varphi \rangle \, dt = 0 \, ,
\end{align}
for it can be easily checked that the limit can be passed in the resulting nonlinear term\footnote{In fact, the limit in equation \eqref{lim_R^n} holds for \emph{any} $u \in L^2_tL^2_x$, regardless whether it is a solution to the Navier-Stokes equations or not. So, the fact that the term involving the cumulant, $- \nabla \cdot R^n$, appears on the RHS of the coarse-grained equations is an additional property of the evolution equations that can potentially be exploited beyond the limit \eqref{lim_R^n}. See \cite{duchon2000inertial}. }.  However, since it is unknown whether the sequence $\nabla \cdot \overline{R}^{ \, n} $ remains uniformly bounded for instance in the natural space $L^2(0,T; V')$, condition \eqref{lim_R^n} alone cannot guarantee the weak convergence of $\nabla \cdot \overline{R}^{ \, n}$ in this space.  
\smallskip
\\
\indent As will be shown, when the weak limit of the sequence of minimizers $v^n_*$ of functional \eqref{functional-intro}, which is guaranteed to exist, does not identically vanish, then the limit  in equation \eqref{lim_R^n} holds for $\varphi \in L^2(0,T;V)$. Key in the proof is replacing, for each $n \in \mathbb{N}$, the sum of the two dissipative terms in the coarse-grained Navier-Stokes equations \eqref{NSCG-intro}, that is $- \nu \Delta \overline{u}^{ \, n}   + \nabla \cdot \overline{R}^{ \, n} $, with a single term possessing improved overall convergence properties. This term is obtained from the Euler-Lagrange equations associated with the minimization problem \eqref{min-prob-intro}:
\begin{equation}    \label{poisson-intro}
(1-2\lambda_n) \Delta v^n_*(x,t) = -\nu \Delta  \overline{u}^{ \, n} + \nabla \cdot \overline{R}^{ \, n} (x,t)\, ,
\end{equation}
understood to hold in the weak sense in space-time. The sequence of real numbers $\lambda_n$ represents the Lagrange multipliers that result when taking into account the constraint, \emph{i.e.} the surface of the enstrophy ball. As a consequence of the system of equations \eqref{poisson-intro} and under the condition $v_* \neq 0$, we obtain $ (1-2\lambda_n) \nabla v^n_* \rightharpoonup \nu \nabla u$ in $L^2_t L^2_x$ as $n \to \infty$, which finally produces the weak convergence of $\nabla \cdot \overline{R}^{ \, n} $ along a subsequence. See Proposition \ref{prop-w-conv}.

\begin{remark} \label{Other variational approaches}
 A variational resolution of weak solutions to the Navier-Stokes was given in \cite{ghoussoub2007antisymmetric, ghoussoub2009anti} by using the theory of anti-self-dual Lagrangians. The evolution equations are derived not through the standard procedure (\emph{i.e.} as the Euler-Lagrange equations), but rather by recognizing that the minimizer is also a zero of the Lagrangian itself. These techniques have also been recognized to be applicable to other stationary and evolutionary problems. In \cite{ortiz2018variational}, a more standard variational approach was developed for deriving Leray-Hopf solutions using a minimization principle. The functional considered is a particular instance of a broader class of functionals known as Weighted Inertia Dissipation Energy (WIDE) functionals. For the problem of incompressible viscous flow, the WIDE functional becomes 
\begin{equation} \label{Wide}
    I^\epsilon(w) = \int_0^\infty \int_\Omega e^{-\frac{t}{\epsilon}} \left( \dfrac{1}{2} | \partial_t w + (w \cdot \nabla) w |^2 + \dfrac{\nu}{2 \epsilon} | \nabla w|^2 + \dfrac{\sigma}{2} |  (w \cdot \nabla) w |^2 \right)\, dx dt \, .
\end{equation}
 The functional is the sum of the inertia and dissipation of the fluid, in addition to a stabilizing term (involving the parameter $\sigma$), which is needed to guarantee the convergence of minimizer $w^\epsilon$ in the limit. By formally computing the Euler-Lagrange equations for $I^\epsilon$ one obtains the incompressible Navier-Stokes equations in $w^\epsilon$ at the \emph{zeroth order} in the parameter $\epsilon$. As $\epsilon \to 0$, the sequence of minimizers $w^\epsilon$ is proven to converge weakly to a Leray-Hopf solution, and the higher order terms vanish in this limit. It is worth remarking that this method, which dates back to the pioneering work of De Giorgi \cite{de1996conjectures}, has a wide array of applications to nonlinear evolution equations of parabolic and hyperbolic types\footnote{In the absence of a rigorous variational principle for hyperbolic PDEs, De Giorgi, while studying a family of semi-linear wave equations, conjectured an elliptic-type regularization method that can be derived from a minimization principle. The desired solutions can then be obtained in the limit of the sequence of minimizers. The conjecture was resolved in \cite{serra2012nonlinear}.}. See \cite{ortiz2018variational} and references therein for more details. 
\medskip
\\
\indent By way of contrast to these approaches, we do not ask the minimizer $v^n_*$ of functional $\mathcal{K}^n$ to solve the Navier-Stokes equations at any order. In this sense, our Euler-Lagrange equations do not represent an approximation scheme to the original equation -- a role customarily assigned to these equations in the context of finding variational solutions to the equations of motion. We, on the other hand, only approximate the dissipative terms by the variational procedure, recognizing that the limit can already be passed in the inertial term. This approximation is, in fact, all that is needed to obtain the desired energy equalities (at least in the $v_* \neq 0$ case) since it implies that the weak limit of $\nabla \cdot R^n$ vanishes on the entire space and not only over a dense subset, as highlighted above.
\end{remark} 

\begin{remark} 
Before concluding this section, we note that the variational method outlined above is applicable, more broadly, to PDEs of the form:
\begin{equation} \label{PDE-gen}
    \partial_t w + \nabla_x \cdot [ \mathbf{F}(w) ] - \Delta w = 0 \, ,
\end{equation}
on $\mathbb{T}^{k_1} \times [0,T]$, and $w \in \mathbb{R}^{k_2}$ and the nonlinear function $\mathbf{F} \in \mathbb{R}^{k_1 \times k_2}$, with $k_1, k_2 \in \mathbb{N}^+$, satisfies $\langle \mathbf{F}(w), \nabla w \rangle = 0$. Now suppose the existence of  a global-in-time weak solution $w$ to equations \eqref{PDE-gen}, for instance in $  C_{\mathrm{weak}} (0, T; H) \cap L^2(0, T; V) $, satisfying an energy inequality similar to \eqref{NS_EN_ineq}:
$$ \dfrac{1}{2}\lVert w(t) \rVert^2 + \int_0^t \lVert \nabla w \rVert^2 ds \leq \dfrac{1}{2}\lVert w(0) \rVert^2\, , $$
for all $t \geq 0$. Let $\overline{w}^{ \, n}$ represent a suitable mollification of the solution $w$, then the mollified field solves the following coarse-grained equations:
\begin{equation*} \label{PDE-gen-moll}
    \partial_t \overline{w}^{ \, n} + \nabla_x \cdot [ \mathbf{F}(\overline{w}^{ \, n}) ] - \Delta \overline{w}^{ \, n} + \nabla_x \cdot [ \overline{\mathbf{Q}(w)}^{ \, n}  ] = 0 \, ,
\end{equation*}
where $ \overline{\mathbf{Q}}^{ \, n}$ is defined analogously to $\overline{R}^{ \, n}$ in equations \eqref{NSCG-intro}. Therefore, as long as $\mathbf{F}(\overline{w}^{ \, n}) \to \mathbf{F}(w)$ as $n \to \infty$ weakly on some dense subset of the solution space and the analogous limit $w_* \neq 0$, we can apply the variational procedure to obtain the global energy equality for the PDE:
$$ \dfrac{1}{2}\lVert w(t) \rVert^2 + \int_0^t \lVert \nabla w \rVert^2 ds = \dfrac{1}{2}\lVert w(0) \rVert^2\, , $$
for all $t \geq 0$.
\end{remark}

\section{Variational analysis} \label{Variational analysis}

In this section, we prove the formulation of Onsager's principle of the least dissipation of energy as outlined in the Introduction (Proposition \ref{prop2.1} below). As such, the central object of study in this section is the following functional

\begin{equation*} 
\mathcal{K}(v) = \int_0^T \frac{1}{2} \lVert \nabla v(t) \rVert^2 - \langle J(t), \nabla v (t) \rangle\, dt\, .
\end{equation*}
We  suppose that we are given a non-trivial sequence of fluxes $J^n: \mathbb{T}^3 \times [0,T] \to \mathbb{R}^{3 \times 3}$  belonging to $L^2_tL^2_x$. Such a sequence will be given in terms of the mollification of  weak solutions of the Navier-Stokes equations, which we consider next.
\subsection{Mollification} \label{mollification}
We begin by defining the notion of weak solutions to the Navier-Stokes equations in the space-periodic setting as follows: let $H $ and $V$ be the closures of the space of periodic, smooth, divergence free, zero mean vector-valued functions in $L^2(\mathbb{T}^3)$ and $H^1(\mathbb{T}^3)$, respectively. Denote $\langle \cdot, \cdot \rangle$ the inner product in $L^2$ space, and $V'$ the natural dual of $V$ with $V \subset H \subset V'$. Then:

\begin{definition} \label{def}\emph{ We say $u\in L^\infty(0,T; H) \cap L^2(0,T; V) $ is a weak solution of \eqref{NS} satisfying the initial condition $u_0 \in H$ if $u(x, t)$ satisfies }
\begin{equation} \label{def-NS}
    -\int_0^T \langle u, \partial_t \varphi \rangle \, dt + \int_0^T \langle (u \cdot \nabla) u, \varphi  \rangle \, dt + \nu \int_0^T \langle \nabla u, \nabla \varphi  \rangle \, dt = \langle u_0, \varphi(0) \rangle \, , 
\end{equation}
\emph{for all divergence free test functions $\varphi \in C^\infty_c (\mathbb{T}^3 \times [0, T))$, and $u(0) = u_0$ weakly in $V'$.}
\end{definition} 
Although we have not required weak solutions to be $L^2$-weakly continuous in-time, this property follows from the definition; see \cite{robinson2016three}. As was mentioned earlier, if, additionally, a weak solution satisfies the global energy inequality \eqref{NS_EN_ineq} for all $t \geq 0$ then it is called a Leray-Hopf solution. In this section and for the rest of the paper we will reserve the designation $u = u(x,t)$ to denote a Leray-Hopf solution on $\mathbb{T}^3 \times [0,T]$.   
\smallskip
\\
\indent Let $\eta^n = \frac{1}{\delta^3_n} \eta \left( \frac{x}{\delta_n} \right) $ represents a family of spatial mollifiers in $\mathbb{R}^3$ introduced earlier in Section 1, where $\eta \in C^\infty$  be a standard smooth nonnegative radial function (kernel) supported over some compact subset of $\mathbb{T}^3$. Then, on $\mathbb{T}^3 \times [0,T] $ the  \emph{filtered} or \emph{coarse-grained} velocity $\overline{u}^{ \, n} = \eta^n * u$ satisfies
    \begin{equation}\label{NSCG-sec 2}
  \left\{
  \begin{array}{ll}
  \vspace*{0.1 in}
   \partial_t \overline{u}^{ \, n} + (\overline{u}^{ \, n} \cdot \nabla) \overline{u}^{ \, n}+ \nabla \overline{p}^{ \, n} - \nu \Delta \overline{u}^{ \, n}  + \nabla \cdot \overline{R}^{ \, n} = 0 \, , \\
    \nabla \cdot \overline{u}^{ \, n} = 0\, ,
    \end{array}
     \right.
\end{equation}
where the \emph{cumulant} associated with the coarse-graining is defined by $\overline{R}^{ \, n}(x,t) = (\overline{u \otimes u})^{n} - \overline{u}^{ \, n} \otimes \overline{u}^{ \, n}$, and the scalar \emph{pressure} field $\overline{p}^{\, n}(x,t)$ is smooth in the spatial variable and satisfies
$ - \Delta \overline{p}^{\, n}( \cdot, t) = (\nabla \otimes \nabla) : (\overline{u \otimes u})^{n}(\cdot, t) $
for almost every $t \in (0,T]$. We will also refer to $\overline{R}^{ \, n}(x,t)$ as the \emph{Reynolds stress} for reasons that will become clear in Section \ref{discussion}. 
\begin{remark} \label{C_w remark}
By double mollification in space-time, we can pass the limit on the LHS of the resolved energy equation \eqref{NSCG_Glob_EN}. Let $\overline{u}^{ \, n}(x,t) = \zeta^n *_{x,t} u $. The space-time mollification is defined as $\zeta^n *_{x,t} \, f = j^n *_t (\eta^n *_x f)$ for $f \in L^1_\mathrm{loc} (\mathbb{R}^3 \times \mathbb{R})$, where $j(t)$ is a another mollifier: a positive, even, smooth function with  unit mass and compact support on $\mathbb{R}$; and functions $j^n$ are defined analogously to $\eta^n$. The mollification in-time is given by
$$ \overline{h}^{ \, n}(t) = \int_0^T j^n(t-s) h(s)\, ds\, , $$
for $h \in L^1_\mathrm{loc}(\mathbb{R})$. Using the fact that solution $u$ is weak-in-time-$L^2$ continuous, we obtain \cite{galdi2000introduction, berselli2023energy}
$$ \dfrac{1}{2} \lVert  u(t) \rVert^2 - \dfrac{1}{2} \lVert u_0 \rVert^2 = \lim_{n \to \infty} \int_0^t \langle \partial_t u, \zeta^n *_{x,t} (\zeta^n *_{x,t} u) \rangle \, ds = \lim_{n \to \infty} \left( \dfrac{1}{2} \lVert \overline{u}^{ \, n}(t) \rVert^2 - \dfrac{1}{2} \lVert \overline{u}^{ \, n}_0 \rVert^2 \right) \, ,$$
for all $t \geq 0.$
\end{remark}

\subsection{Minimization problem} \label{Minimization problem}

In this subsection our goal is to produce a solution the minimization problem \ref{min-prob-intro}. Since $\xi \to K(x,t,\xi) \doteq \frac{\nu}{2} \xi^2 - J^n(x,t) : \xi $ is strictly convex for almost every $(x,t) \in \mathbb{T}^3 \times [0,T]$, the density function $K(x,t,\xi)$ satisfies the standard coercivity and growth estimates, and due to the quadratic nature of the constraint; we can use the direct method of the calculus of variations over the enstrophy ball $\mathcal{B} = \lbrace v \in L^2(0,T; V) :  \int_0^T \lVert \nabla v \rVert^2dt \leq \frac{1}{2} \lVert u_0 \rVert^2 \rbrace$ to obtain 
\begin{proposition} \label{prop2.1} 
Fix some $n \in \mathbb{N}^+$. Let 
\begin{align} \label{functional-sec3}
    \mathcal{K}^n(v) = \int_0^T \frac{1}{2} \lVert \nabla v(t) \rVert^2 - \langle J^n(t), \nabla v (t) \rangle\, dt \, .
\end{align}
Then
\begin{itemize}
    \item[1.] There exists a unique function $v^n_* \in \mathcal{B}$ such that $\mathcal{K}^n(v^n_*) = \inf \lbrace \mathcal{K}^n(v) : v \in \mathcal{B} \rbrace$. 
    \item[2.] There exists a real number $\lambda = \lambda_n$ such that the constrained Euler-Lagrange equations are satisfied:
    \begin{equation} \label{poisson-weak}
\int_0^T \langle  \nabla v_*^n (t), \nabla \varphi(t) \rangle\, dt = \int_0^T \langle  J^n(t), \nabla \varphi(t) \rangle\, dt + 2 \lambda_n \int_0^T \langle \nabla v_*^n, \nabla \varphi(t) \rangle  \,dt \, ,
    \end{equation}
    for all divergence free $\varphi \in C^\infty (0, T; C^\infty_c)$.
\end{itemize}
\end{proposition}
We include the proof of the above proposition in Appendix \ref{Appendix B} for the readers' convenience.
\begin{remark} \label{remark_multiplier}
 In general, if $\lambda_n \neq 0 $ (implying $v^n_* \in \partial\mathcal{B}$) then we have 
\begin{equation} \label{lambda_eq}
(1-2\lambda_n) = \dfrac{ \int_0^T \langle  J^n(t), \nabla \psi(t) \rangle\, dt}{\int_0^T \langle  \nabla v_*^n (t), \nabla \psi(t) \rangle\, dt} \, ,
\end{equation}
for some (non-trivial) smooth, divergence free test function $\psi$; see Remark \ref{remark-def_lagrange_multiplier}.
\end{remark}

Since the sequence $J^n(x,t)$, defined by equation \eqref{sigma}, is smooth in the spatial variables, the following result immediately follows
\begin{corollary} \label{cor}
The minimizer $v^n_*$ satisfies the following Poisson equations point-wise in $x\in \mathbb{T}^3$ and distributionally in $t \in [0,T]$:
\begin{align} \label{poisson-sec3}
- (1-2\lambda_n) \Delta v^n_*(x,t) = - \nu \Delta \overline{u}^{\, n}(x,t) +  \nabla \cdot \overline{R}^{ \, n}(x,t)  \, .
\end{align}
\end{corollary}
The system of Poisson equations \eqref{poisson-sec3} together with the boundedness of the sequence of minimizers in $L^2(0,T;V)$ will be key in passing the limit in the coarse-grained energy equations, as discussed next.

\section{Passage to the limit} \label{passage}

The approach advanced in the previous section has furnished a sequence (of minimizers) $v^n_*$ uniformly bounded in $L^2(0,T; V) $. However, in view of the structure of the Poisson equation \eqref{poisson-weak} and equation \eqref{lambda_eq} containing the Lagrange multipliers, we exclude the case where the weak limit, which we denote by $v_*$, vanishes. See Section \ref{discussion} for further discussion on this case.
\smallskip
\\
\indent We recall that establishing the existence of a subsequence satisfying $\lim_{n \to \infty} \langle  \overline{R}^{\, n}, \nabla \varphi \rangle_{L^2_t L^2_x} = 0$ for all $ \varphi \in L^2(0,T; V)$ amounts to the energy equalities \eqref{global-en-eq-intro} and \eqref{local-en-eq-intro}. This follows directly from the next result, which is a simple consequence of Proposition \eqref{prop2.1} and our assumption. We include the proof for the sake of completeness. 

\begin{proposition} \label{prop-w-conv}
Assume $v_* \neq 0$. If $J^n(x,t)$ is given by equations \eqref{sigma}, then $ (1-2\lambda_n) v^n_* \rightharpoonup \nu u $ in  $ L^2(0,T; V)$, as $n \to \infty$ (along a subsequence).
\end{proposition}
{\bf Proof}.
\medskip
\\
\indent Consider $\varphi \in L^2(0,T; V)$ and a sequence of divergence free functions $\varphi_m \in C^\infty_c (\mathbb{T}^3 \times [0, T))$ such that $\nabla \varphi_m \to \nabla \varphi  $ strongly in $L^2_tL^2_x$ (with $n > m$). Then
\begin{align*}
    \int_0^T \langle (1-2\lambda_n)\nabla v^n_* - \nu \nabla u, \nabla \varphi \rangle \, dt &\leq \left( |1-2\lambda_n| \lVert \nabla v^n_*\rVert_{L^2_tL^2_x} + \nu \lVert \nabla u \rVert_{L^2_tL^2_x} \right)  \lVert \nabla \varphi_m  - \nabla \varphi \rVert_{L^2_tL^2_x}\\
 &\qquad \qquad \qquad \qquad \qquad \qquad+ \int_0^T \langle (1-2\lambda_n) \nabla v^n_* - \nu \nabla u, \nabla \varphi_m  \rangle \, dt \, ,\\
   &\leq C  \lVert \nabla \varphi_m  - \nabla \varphi \rVert_{L^2_tL^2_x} +  \int_0^T \langle (1-2\lambda_n) \nabla v^n_* - \nu \nabla u, \nabla \varphi_m  \rangle \, dt \, ,
\end{align*}
where constant $C$ is independent of $n$ and $m$. The last inequality follows since $v^n_*$ belongs to $\mathcal{B}$ for every $n \geq 1$, and $ (1 - 2\lambda_n)$ is convergent hence bounded (along a subsequence) by equation \eqref{lambda_eq} (note that by equation \eqref{lim_R^n}, the numerator in the RHS of equation \eqref{lambda_eq} converges for smooth, divergence free test functions). We conclude by sending $n, m \to \infty$ and invoking equations \eqref{poisson-weak} and \eqref{lim_R^n}.

\begin{remark} \label{remark-direct method}
Following Corollary \ref{cor}, the coarse-grained equations can be expressed as
\begin{equation*}
    \partial_t \overline{u}^{\, n} + \nabla \cdot [\overline{u}^{\, n} \otimes \overline{u}^{\, n}] + \nabla \overline{p}^{\, n}  -  (1 - 2\lambda_n) \Delta v^n_* = 0  \, ,
\end{equation*}
with energy equation 
\begin{equation} \label{en_balance_min}
    \frac{1}{2} \lVert \overline{u}^{\, n}(T) \rVert^2 - \frac{1}{2} \lVert \overline{u}^{\, n}(0) \rVert^2   = -  (1 - 2\lambda_n)  \int_0^T \langle  \nabla v^n_*, \nabla \overline{u}^{\, n} \rangle \, dt \, ,
\end{equation}
for all $T>0$. We can directly prove the global energy equality \eqref{global-en-eq-intro} using the energy balance \eqref{en_balance_min} without reference to equation \eqref{NSCG_LOC_EN}, since $ (1 - 2\lambda_n) \nabla v^n_* \rightharpoonup \nu \nabla u$ in $L^2_tL^2_x$ as $n \to \infty$, and the limit can be passed on the LHS of equation \eqref{en_balance_min} (see Remark \ref{C_w remark}). 
\end{remark}

\begin{remark} \label{remark_compare}
We conclude by making explicit a few additional limiting properties of the Reynolds stress $\overline{R}^{\, n}$. For simplicity, we set $\nu = 1$. Since both $v^n_*$ and $- v^n_*$ belong to the admissible set $\mathcal{B}$, the comparison $\mathcal{K}^n(v^n_*) < \mathcal{K}^n(- v^n_*)$ yields
\begin{align} \label{R^n_property1}
   \lim_{n \to \infty} \int_0^T \langle \overline{R}^{\, n}, \nabla v^n_* \rangle \, dt \leq \int_0^T \langle \nabla u, \nabla v_* \rangle \, dt \, ,
\end{align}
along a subsequence. By invoking $\mathcal{K}^n(v^n_*) < \mathcal{K}^n(- v^n_*)$ again, we get 
$$ \int_0^T  \langle \nu \nabla \overline{u}^{\, n}  - \overline{R}^{\, n}, \nabla v^n_* \rangle \, dt > 0\, \quad \forall n , $$
which, when combined with the constrained Euler-Lagrange equations, implies
$$ \lim_{n \to \infty} (1- 2\lambda_n) \geq 0\, , $$ 
along a subsequence. In this case, a sharper estimate can be derived by comparing $v^n_*$ with $\overline{u}^{\, n}$ arising from Leray solutions, which even yields a definite sign for the limiting value in \eqref{R^n_property1}.
\begin{align*}
\dfrac{1}{4} \lVert u_0  \rVert^2 + \int_0^T \langle \overline{R}^{\, n}, \nabla v^n_* \rangle \, dt & < \int_0^T \langle  \nabla v^n_*, \nabla \overline{u}^{\, n} \rangle \, dt - \frac{1}{2} \int_0^T \lVert \nabla \overline{u}^{\, n} \rVert^2 dt + \int_0^T \langle \overline{R}^{\, n}, \nabla \overline{u}^{\, n}  \rangle \, dt \, , \\
\limsup_{n}  \int_0^T \langle \overline{R}^{\, n}, \nabla v^n_* \rangle \, dt & \leq - \frac{1}{2} \int_0^T \lVert \nabla u \rVert^2 dt - \int_0^T \int_{\mathbb{T}^3} D(u)\, dx dt -\dfrac{1}{4} \lVert u_0  \rVert^2 \, .
\end{align*} 
Finally, drawing on Corollary \ref{cor} and Proposition \ref{prop-w-conv}, when $v_* \neq 0$ the following convergence properties hold:
$$ \limsup_n\, \lVert \nabla \cdot \overline{R}^{\, n} \rVert_{L^2(0,T;V')} < \infty\, , \qquad \lim_{n \to \infty} \int_0^T \langle \overline{R}^{\, n}, \nabla u \rangle\, dt =0 \, . $$  
\end{remark}

\section{Discussion} \label{discussion}

A seemingly simple and innocuous criterion (\emph{i.e.} $v_* \neq 0$ in Theorem \ref{Th1}) for energy conservation can be understood and interpreted within the frameworks developed for the study of turbulent flows from at least two distinct, yet related, perspectives. First, the Poisson equations \eqref{poisson-sec3} can be regarded as prescribing a constitutive law -- or, more precisely, a modeling assumption -- for the cumulant or the Reynolds stress $\overline{R}^{\, n} (x,t)$. From a physics point of view, this is indeed one of the primary implications of the principle of the least dissipation of energy as formulated by Onsager (see Section \ref{main result} and Appendix \ref{Appendix A}). In essence, equations \eqref{poisson-sec3} suggest that the Reynolds stress may be represented as follows (up to a divergence):

\begin{equation} \label{reynolds-modeling}
\overline{R}^{\, n} (x,t) = 2 (1-2\lambda_n) \, \mathrm{sym}( \nabla v^n_* )(x,t) \,.
\end{equation}
In other words, within the framework of irreversible thermodynamics outlined in Appendix \ref{Appendix A}, the contribution to internal entropy production of the system arises predominately from the thermodynamic fluxes associated with velocity fluctuations, namely $|\mathbf{J}| \approx |\mathbf{J}_2|$. We note that equations \eqref{reynolds-modeling} are analogous to the (deviatoric) \emph{Boussinesq eddy-viscosity approximation}, where the mean strain rate is replaced by a suitable approximation, and the prefactor involving the Lagrange multiplier plays the role of an effective eddy viscosity (see \cite{wilcox1998turbulence} for a modern exposition on turbulence modeling). The limit $v^n_* \rightharpoonup 0$ in $L^2(0,T;V)$ therefore corresponds to larger and larger values of "eddy viscosity", as can be seen by equation \eqref{lambda_eq}. This, in turn, causes the $L^2_tL^2_x$ norm of the Reynolds stress to grow unstable as can be seen from equation \eqref{reynolds-modeling}: 
$$ | 1-2 \lambda_n| \lVert u_0 \rVert \leq \sqrt{2 \int_0^T \lVert \nabla \overline{R}^{\, n}\rVert^2 dt} \, .$$
 To the extent that the Boussinesq approximation remains valid, these observations imply that the flow is being modeled in a highly turbulent regime with eddy-driven transport dominating molecular diffusion. This does not necessarily imply that energy conservation \eqref{global-en-eq-intro} is violated, but it does point to a potential deficiency in our variational approach: the sequence $v^n_*$ does not carry enough information about the energy or the energy dissipation rate of the flow, an interpretation which we can associate with the indeterminacy observed in the case $v_* =0$. 
\medskip
\\
\indent To further illustrate this point, we present one possible way to improve the properties of the sequence of minimizers associated with the functional \eqref{functional-sec3}. The idea is to introduce an additional degree of freedom into the problem by prescribing the total energy profile. This method was considered in \cite{CHEN2024134274}, as mentioned in Section \ref{introduction}, however in the context of a probabilistic study of the Navier-Stokes equations. In particular, it was conjectured the existence of a probability measure $\eta_{e,t}$  supported over the fixed-in-time \emph{energy-enstrophy surface}:
\begin{equation} \label{surface}
    G_e(t) = \lbrace v \in L^\infty_t L^2_x \cap L^2_t H^1_x \, : \, \dfrac{1}{2} \lVert v(t) \rVert^2 + \nu \int_0^t \lVert \nabla v(s) \rVert ^2 ds  = e(t) \rbrace \, ,
\end{equation}
where $e(t)$ is a given positive scalar -- to be understood as the prescribed energy profile\footnote{The significance of  measure $\eta_{e,t}$, termed the physical measure, stems from the fact that it maximizes, for every $t \in [0,T]$, the relative entropy functional $ \mathcal{S}(f_{e,t}) =  - \displaystyle\int f_{e,t}(v) \log (f_{e,t}(v)) \, d\eta_{e,t}(v) $ where $f_{e,t}(v) := \frac{d \mu_{e,t}(v)}{d \eta_{e,t}(v)} \in L^1(\eta_{e,t})$. This flexibility in determining the maximum entropy measure was used to draw connections to the Kolmogorov theory of turbulence.}.
\medskip
\\
\indent In the deterministic setting, a connection can be made with the principle of the least dissipation of energy by considering a variant of the energy-enstrophy surface as the admissible set to the minimization problem. To this end, define 
\begin{equation} \label{EE-surface}
\mathcal{G}_e(t) = \lbrace v\in  L^\infty_t L^2_x \cap L^2_t H^1_x : \,  \lVert v \rVert_t^2 \leq e(t) \rbrace \, ,
\end{equation}
where $\lVert v \rVert_t^2 = \frac{1}{2} \lVert v(t) \rVert^2 + \nu \int_0^t \lVert \nabla v(s) \rVert ^2 ds$ and function $e(t)$ is prescribed in advance. One can naturally include the initial condition $u_0$ into the formalism by defining the following admissible set, consisting of functions with higher regularity than those in $\mathcal{B}$:
\begin{align}
   \mathcal{B}' = \lbrace v \in \mathcal{G}_e(t) : v(x,0) = u_0, \, \, \mathrm{for \, \, a. \, a.} \, \, t \in [0,T] \rbrace \, ,
\end{align}
where the initial condition is satisfied in the sense of $L^2$-weak limits (in time). It might be also desirable to replace $u_0$ with $u_0^m$ -- some suitable approximation to the initial condition. In any case, we require that function $e(t)$ is appropriately chosen so that $\mathcal{B}'$ is not the empty set.  
\smallskip
\\
Following the proof of Proposition \eqref{prop2.1}, we can conclude the existence of a unique minimizer $v^{n,e}_{*}(t) \in \mathcal{G}_e(t)$ to the functional \eqref{functional-sec3} for almost all $t$ (note that relation \eqref{cauchy-seq} implies strong convergence in $L^2(0,T; H)$ and therefore $\frac{1}{2} \lVert v_{m}^{n,e} (t) \rVert^2 \to \frac{1}{2}  \lVert v_*^{n,e} (t) \rVert^2$ along a subsequence for almost every $t$.).
\smallskip
\\
We now investigate the $L^2$-in-time weak limit of $v^{n,e}_{*}$ at $t=0$, that is does $\lim_{t\to0} \langle v^{n,e}_{*}(t) , w \rangle =\langle u_0 , w \rangle $ hold for every $w \in L^2(\mathbb{T}^3)$? Let  $v^{n,e}_{m} \in \mathcal{B}'$ denote a minimizing sequence (in $m$) for the functional $\mathcal{K}^n$. Then, for $t$ close to zero and due to the convexity of $\mathcal{G}_e(t)$, we have:
\begin{align*}
| \langle v^{n,e}_{*}(t) - u_0 , w \rangle | &\leq | \langle v^{n,e}_{*}(t) - v^{n,e}_{m}(t) , w \rangle | + | \langle v^{n,e}_{m}(t) - u_0 , w \rangle | \\
& \lesssim \left \lVert \dfrac{v_{*}^{n,e}(t) - v_m^{n,e}(t)}{2} \right \rVert^2 + O(t^\alpha) \\
& \lesssim e(t) + O(t^\alpha), \label{est_e}\numberthis 
\end{align*}
for some $\alpha \geq 1$. Therefore, the continuity of the $v^{n,e}_*$ near zero is controlled by the prescribed energy $e(t)$, a feature that arises from reformulating the variational principle over $\mathcal{B}'$. In fact, one can also readily derive an estimate analogous to \eqref{est_e}  away from $t=0$. It follows, the $L^2$-weak continuity of $v^{n,e}_*(x,t)$ at any $t \in [0,T]$ is controlled by the energy $e(t)$, not only in a neighborhood of $t=0$. This simple example illustrates a possible interplay between the regularity -- and potentially the convergence -- of the sequence of minimizers $ v_*^{n,e}$ and the properties of $e(t)$. In general, prescribing this energy is expected to be closely linked to controlling the higher Fourier modes of solutions to the Navier-Stokes equations and to guaranteeing their convergence, in the inviscid limit, to solutions of the Euler equations. The first rigorous example of this strategy was presented in \cite{chen2012kolmogorov}. In the same vein, when $\lambda_n \neq 0$ we obtain an explicit expression for the minimizer $v^{n, e}_*$ in terms of the prescribed energy: $\int_0^t 2 \, |\mathrm{sym}(\nabla v^{n,e}_*(x,s)) |^2 \, ds = e_n(t)$. Yet, we expect that merely specifying $e_n(t)$ in terms of the initial energy profile does not produce a velocity that can adequately capture certain essential features of the dynamics at finer scales\footnote{In fact, in 1953 Onsager (together with S. Machlup) \cite{onsager1953fluctuations} extended his original theory in response to a related challenge. While the reciprocal relations predict the average behavior of a system -- for example, the mean heat flux or mean diffusion current -- the underlying dynamics of the fluctuations that give rise to the averages were not taken into account. The extended theory introduced an explicit random Gaussian noise into the formalism, allowing for the calculation of the probability paths of a given state variable.}. It may well prove beyond reach to determine an appropriate sequence of prescribed energies in space-time coordinates capable of overcoming this difficulty. An alternative way to compensate for the loss of regularity observed would be to extend Onsager's variational principle to incorporate higher-order terms. From a physics perspective, this approach is closely connected to calculating the rate(s) of entropy production in systems far from equilibrium (\emph{e.g.}, in flux space), an aspect of non-equilibrium thermodynamics still not well-established even if restricted to isolated systems (see \cite{martyushev2014restrictions}, for instance,  for a critical review on the subject).

\bigskip

\begin{appendices}

\section{Principle of least dissipation of energy} \label{Appendix A}

In this Appendix, we provide a short account of the principle of least dissipation of energy -- first introduced by L. Onsager in his seminal work in linear irreversible thermodynamics \cite{onsager1931reciprocal1, onsager1931reciprocal2} -- within the framework of classical irreversible thermodynamics. The mathematical formulation of the principle will be confined to incompressible isotropic Newtonian fluids. We refer the reader to \cite{extended} for a general exposition on the theory.

\subsection{Overview of theory}

Throughout this section, we assume for simplicity that the fields under consideration belong to $C^\infty(\overline{\Omega} \times [0,T])$. A more general setting will be considered at the end of the discussion. We begin with a body  $\Omega \subset \mathbb{R}^3$ having smooth boundary and mass density $\rho$. We assume that the physical processes taking place in $\Omega$ are determined by the evolution of $m$ state variables: $\mathbf{a}_1(x,t), \mathbf{a}_2(x,t), \cdots, \mathbf{a}_m(x,t)$ obeying conservation laws of the form

\begin{equation} \label{evolv-eq}
    \rho \dfrac{\mathrm{D} \mathbf{a}}{\mathrm{D}t} = \nabla \cdot \boldsymbol\sigma_\mathbf{a} + \mathbf{G}(x,t) \, ,
\end{equation}
where $\mathbf{a}(x,t) = (\mathbf{a}_1(x,t), \mathbf{a}_2(x,t), \cdots, \mathbf{a}_m(x,t))$, $\mathbf{G}(x,t)$  represent body forces and/or source functions, and $\boldsymbol{\sigma}(x,t)$ represent a combination of the fluxes (stresses) acting on and across $\Omega$; we drop the subscript $\mathbf{a}$ for convenience. The state variables could either be scalar or vector quantities, in which case $\boldsymbol{\sigma}$ correspond to vector or tensor fields, respectively. We write $\boldsymbol\sigma$ as the sum of its (quasi)conservative and dissipative parts:
\begin{align*}
    \boldsymbol\sigma =\boldsymbol\Pi + \mathbf{J} \, ,
\end{align*}
where each competent of the dissipative stresses $\mathbf{J}_k(x,t)$, $k=1, 2, \cdots, m$, represents a generalized thermodynamic flux acting within $\Omega$ and is associated with a distinct irreversible process given by either a single state variable $\mathbf{a}_k$ or a subset of them. 
\medskip
\\
\indent On the other hand, the local balance of entropy in  $\Omega$ reads
\begin{align} \label{entropy_balnce}
    \rho  \dfrac{\mathrm{D}s}{\mathrm{D}t} + \nabla \cdot \mathbf{h}_s = \rho \dot{s}^{(i)} \, ,
\end{align}
where $s(x,t)$ and $\dot{s}^{(i)}(x,t)$ are the specific entropy and internal entropy production of the system, respectively, and $\mathbf{h}_s$ represents the entropy flux across the boundary of $\Omega$.   
\medskip
\\
\indent The dependence of the entropy of the system on the state variables becomes clear once we write the function $\dot{s}^{(i)}$ in the standard form as the product of the thermodynamic fluxes and forces:
\begin{align} \label{entropy_product}
    \rho \dot{s}^{(i)} =\rho \dot{s}^{(i)}(\mathbf{J}, \mathbf{X}) = \sum_k \mathbf{J}_k \mathbf{ X}_k \, .
\end{align}
Thermodynamic forces $\mathbf{X}_k(x,t)$ are typically known quantities associated with the gradients of the state variables $\mathbf{a}_k$, which may represent, for instance, velocity, concentration, chemical potential, temperature, etc; and are not forces in the mechanical sense. In this context, the second law of thermodynamics is extended to assert that $\dot{s}^{(i)}(x,t) \geq 0$ for all $(x,t) \in \Omega \times (0,T)$. This inequality goes beyond the conventional statement of the second law in equilibrium thermodynamics, which considers only the overall increase in entropy between two equilibrium states of an (isolated) system. It is quite remarkable that the bilinear structure of the internal entropy production function, in terms of forces $\mathbf{X}_k$ and fluxes $\mathbf{J}_k$, is preserved across a wide range of physical systems outside of their equilibrium states. 
\medskip
\\
\indent However, while thermodynamic fluxes are usually not known quantities, they cannot be assigned arbitrarily. In 1931, Onsager argued that thermodynamic processes are realized if they  maximize the following functional in flux space for known values of forces $\mathbf{X}_k = \mathbf{X}_k(x,t)$:
\begin{align} \label{max_principle_appendix}
     \underset{\hspace{-3.1cm} \Omega \times [0,T]}{\int \rho \left( \dot{s}^{(i)}(\mathbf{J}, \mathbf{X}) - \Phi (\mathbf{J}) \right)} \, dx dt \, ,
\end{align}
where the \emph{dissipation function} is given by the following quadratic form
\begin{align*}
   \rho \Phi (\mathbf{J}) = \dfrac{1}{2} \sum_{k,n} \mathbf{R}_{kn} \mathbf{J}_k \mathbf{J}_n \, ,
\end{align*}
and $\mathbf{R}$ is a constant positive (generalized) matrix, and can always be made symmetric. The motivation behind the above variational principle formulated by Onsager becomes evident once we (formally) carry-out the calculation of the first variation. Assume that the functional \eqref{max_principle_appendix} attains its maximum at $\mathbf{J}$, then the Euler-Lagrange equations give

\begin{equation*}
\mathbf{X}_k = \sum_n \mathbf{R}_{kn} \mathbf{J}_n \, ,
\end{equation*}
which implies 
\begin{equation} \label{Ons_rel1}
\mathbf{J}_k = \sum_n \mathbf{L}_{kn} \mathbf{X}_n \, ,
\end{equation}
where $\mathbf{R} = \mathbf{L}^{-1}$ provided matrix $\mathbf{R}$ is invertible. If we identify $\mathbf{L}$ with the matrix of phenomenological coefficients, then equations \eqref{Ons_rel1} together with the resulting symmetry of $\mathbf{L}$ are Onsager's reciprocal relations for linear non-equilibrium thermodynamics. Equations \eqref{Ons_rel1} are examples of \emph{constitutive relations}, which, are essential for closing -- and ultimately solving -- the system of evolution equations \eqref{evolv-eq}, once the initial and/or boundary conditions are specified.
\medskip
\\
\indent  In general, since both the fluxes and the forces depend on the state variables $\mathbf{a}_k$'s and their derivatives, it is natural to expect a functional relation between them. In the equilibrium state of the system, no irreversible processes take place by definition, and so $\mathbf{J}_k = \mathbf{J}_k^0 =0$ for all $k$. The linear equations above can thus be naturally interpreted as a first-order expansion of the fluxes around $\mathbf{J}_k^0$:
$$ \mathbf{J}_k = \sum_n \left( \dfrac{\partial \mathbf{J}_k}{\partial \mathbf{X}_n}\right)_{|0 } \hspace{-2.7mm} \mathbf{X}_n + \mathrm{o} (|\mathbf{J}|) \, ,$$
where the matrix of phenomenological coefficients is given by $ \mathbf{L}_{kn} = \left( \frac{\partial \mathbf{J}_k}{\partial \mathbf{X}_n}\right)_{|0 }$, and the evaluation corresponds to the equilibrium state. All the principal minors of the matrix are required to be non-negative to ensure consistency with the second law\footnote{Onsager derived the reciprocal relations based on considerations arising from statistical mechanics;  namely the property of microscopic \emph{time reversal invariance}, which is a statement about the invariance of the evolution equations governing the particles with respect to time reversal. In other words, the principle implies that trajectory of particles are reversed when tracing them backward through time. It was from this microscopic principle that Onsager obtained the macroscopic laws.}.

\subsubsection*{Mathematical formulation for dissipative fluids}

We are now in a position to formulate the functional \eqref{functional-intro}, which is specific to a Navier-Stokes fluid. Consider a dissipative isothermal fluid with velocity field $v$ that vanishes on the boundary $\partial \Omega$, and the fluid temperature normalized to 1. We interpret velocity gradients as thermodynamic forces, so the internal entropy production function reads
\begin{align*}
   \rho s^{(i)}(x,t) = \mathbf{J}: \mathrm{sym} (\nabla v) \, ,
\end{align*}
where $\mathbf{J}$ represents the (symmetric) dissipative stresses acting within the fluid\footnote{It is worth noting, as pointed out in \cite{extended}, that the roles of fluxes and forces can be interchanged without affecting the final outcome. For example, see \cite{ziegler2012introduction}, where this reversal is employed in the context of dissipative thermal fluids. }. We assume that two known or given distinct thermodynamic processes $\mathbf{J}_1$ and $\mathbf{J}_2$ contribute to entropy production. The first is conventional dissipation through viscous effects, which for an incompressible isotropic Newtonian fluid is taken to be $\mathbf{J}_1 = 2 \nu \, \mathrm{sym}(\nabla u)$ where $\nu$ is the constant kinematic viscosity and $u$ represents a solution \emph{a priori} known. The second contribution $\mathbf{J}_2$ captures possible dissipation due to turbulent fluctuations in the velocity field $u$.
\medskip
\\
\indent Let $\mathcal{B}$ be some appropriate admissible set for this setting (\emph{e.g.} a subset of $L^2_tH^1_x$), then the principle of least dissipation of energy \eqref{max_principle_appendix} can be recast into the following conjugate problem 
\begin{align} \label{functionl_final}
     \inf_{v \in \mathcal{B}}  \int_0^T \int_\Omega \left( \dfrac{\nu}{2} | \nabla v |^2 - \mathbf{J} : \nabla v \right) dx dt\, .
\end{align}
(note $\int_\Omega |\nabla v|^2 dx =  2\int_\Omega |\mathrm{sym}(\nabla v)|^2 dx$ on divergence free fields vanishing at the boundary.)

\section{Proof of Proposition \ref{prop2.1}} \label{Appendix B}
Let $\mu^n = \inf \lbrace \mathcal{K}^n(v) : v \in \mathcal{B} \rbrace $. Then, the coercivity of functional $\mathcal{K}^n$ ensures that $\mu^n > -\infty$. We, therefore, select a minimizing sequence $v^n_m \in \mathcal{B}$ (depending on the choice of $u_0$ as well) such that $\mathcal{K}^n(v^n_m) \to \mu^n$, as $m \to \infty$. Again, the standard coercivity estimate gives
\begin{align*}
    \int_0^T \lVert \nabla v^n_m(t) \rVert^2 dt \leq C_1 \, ,
\end{align*}
where the constant $C_1 >0$ is independent of $m$. Therefore, there exists a subsequence of $v_m^n$, which we still call $v^n_m$, and a function $v^n_* \in L^2(0,T; V)$ such that $v^n_m \rightharpoonup v^n_*$ in $L^2(0,T; V)$ as $m \to \infty$.
\smallskip
\\
For the minimizing sequence $v^n_m$, we note the  identity: 
\begin{align*} \label{cauchy-seq}
\frac{1}{8} \lVert v^n_m - v^n_k \rVert^2_{L^2(0;T, V)} &= \frac{1}{2}\mathcal{K}^n(v^n_m) + \frac{1}{2}\mathcal{K}^n(v^n_k) - \mathcal{K}^n \left( \frac{v^n_m + v^n_k}{2} \right) \, , \\
&\leq \frac{1}{2}\mathcal{K}^n(v^n_m) + \frac{1}{2}\mathcal{K}^n(v^n_k) - \mu^n \, , \numberthis 
\end{align*}
for all $k \geq m$. Hence, the minimizing sequence $v^n_m$ is in fact Cauchy in $L^2(0,T; V)$, and so $v^n_m \to v^n_*$ in $L^2(0,T; V)$ as $m \to \infty$. By the Poincare inequality, we additionally have $v^n_m \to v^n_*$ in $L^2(0,T; H)$ as $m \to \infty$, and so the mean is also preserved. Now since $v^n_m \in \mathcal{B} \, \, \forall m \in \mathbb{N}^+$, by the convergence properties of $v_m^n$ we get
\begin{align} \label{en_eq_lim_m}
 \int_0^T \lVert \nabla v^n_*(t) \rVert^2 \, dt \leq \dfrac{\lVert u_0 \rVert^2}{2} \, ,
\end{align}
and so $v^n_* \in \mathcal{B}$. By appealing to the convexity of the of $\mathcal{K}^n$, we conclude that $\mu^n = \mathcal{K}^n(v^n_*).$
\smallskip
\\
Let $v_*$ and $\Tilde{v}_*$ be two distinct minimizers for our problem. Then by the (strict) convexity of $\xi \to K(x,t,\xi)$, we can conclude that $\dfrac{v_* + \Tilde{v}_* }{2}$ is also minimizer, and so $\nabla v_* = \nabla \Tilde{v}_* $ for almost every $(x,t) \in \mathbb{T}^3 \times [0,T]$. Therefore, $v_*$ and $\Tilde{v}_*$ differ by a (time-dependent) constant. But, since $v_*$ and $\Tilde{v}_*$ have zero mean, we conclude that $v_* = \Tilde{v}_* $ a.e. 
\medskip
\\
Write the constraint on the boundary as 
\begin{equation} \label{constraint-proof}
\mathcal{Q}(v) = \int_0^T \lVert \nabla v(t) \rVert^2 \, dt - \dfrac{\lVert u_0 \rVert^2}{2}=0 \, .
\end{equation}
Now define the two scalar functions
\begin{align*}
\Phi (\varepsilon, \tau) = \mathcal{K}^n(v^n_* + \varepsilon \varphi + \tau \psi) \quad \mathrm{and} \quad \Psi (\varepsilon, \tau) = \mathcal{Q}(v^n_* + \varepsilon \varphi + \tau \psi)\, ,
\end{align*}
where $\varphi, \psi \in C^\infty(0,T; C^\infty_c)$, divergence free, and $\psi$ is non-trivial. Due to the convexity of the functional and the constraint, the gradients are well-defined. The constrained Euler-Lagrange equations now follow by the standard Lagrange multiplier theorem applied to functions $\Phi$ and $\Psi$.

\begin{remark} \label{remark-def_lagrange_multiplier}
In relation to Proposition \ref{prop-w-conv}, the test function $\psi(x,t)$ in the proof above can be chosen so that the first variation of $\mathcal{Q}$ with respect to it does not identically vanish for all $n \in \mathbb{N}$ and such that the limit $\int_0^T \langle \nabla v_*, \nabla \psi \rangle \, dt$ is non-zero. Consequently, the Lagrange multiplier
\begin{equation}
\lambda =  \dfrac{\dfrac{d}{d \tau}(\Phi(0, 0))}{\dfrac{d}{d \tau}(\Psi(0, 0))}\, ,
\end{equation}
which depends on $n$ and the choice of the test function $\psi(x,t)$, is well-defined.
\end{remark}

\bigskip
\bigskip
\textbf{Declaration of competing interest}: The author declares that he has no known competing financial interests or personal relationships that could have appeared to
influence the work reported in this paper.
\\
\bigskip
\\
\indent \textbf{Data availability}: No data was used for the research described in the article.
\\
\bigskip
\\
\indent \textbf{Acknowledgments.} The author is grateful to G-Q. Chen, J. Glimm and M. Johnson for helpful discussions and/or comments on this work. The author is also thankful to the referees for their valuable comments and suggestions, which significantly improved the manuscript.
\end{appendices}

\bibliographystyle{siam}
\bibliography{References1_bib}
  
\end{document}